\documentstyle[12pt]{article}

\newtheorem{theorem}{Theorem}[section]

\newtheorem{corollary}{Corollary}[section]
\input{amssym.def}
\input{amssym.tex}

\title{On the Brauer monoid for finite fields}
\author{V.\,V.\,Kirichenko (Kiev, Ukraine),\\
B.\,V.\,Novikov (Kharkov, Ukraine)}
\date{}

\begin{document}
\maketitle

\begin{abstract}

The definition  of the Brauer monoid was given in \cite{hls}. In this
article  it is studied by the notions of modifications \cite{n3} and
0-cohomology \cite{n2}. We investigate the impact of invertible
elements of modifications on the structure of the Brauer monoid,
especially for finite fields.
\end{abstract}

\section{Introduction}

It is well-known that the Brauer group of any finite field is
trivial \cite{dk}. Therefore the so called Brauer monoid proposed in
\cite{hls} is of interest. This monoid generalizes the Brauer group
and isn't trivial for any non-trivial field extension. One can hope
that the studying of its properties will be useful for the investigation
of algebras over finite fields.

The description of the Brauer monoid by modifications an their 0-cohomology
proposed in \cite{n3} is more convenient in our opinion than original one
\cite{hls}. We give in Section 2 this description and a necessary information
about semigroup 0-cohomology too.

Section 3 is devoted to the proof of a theorem, which facilitates
essentially the calculation of a  Brauer monoid for  finite fields by
elimination of the invertible elements.
                         
Finally we note that at stretch of this article relative Brauer monoids
(and relative Brauer groups) are considered only, so the adjective
``relative'' will be omitted.

\section{Preliminary: 0-cohomology and modifications}

Semigroup 0-cohomology  is a specific case of partial cohomologies which
were built in \cite{n2}; it had appeared in the investigation of the
projective representations of semigroups.

Let $S$ be an arbitrary semigroup with a zero. An Abelian group $A$ is
called
{\it a 0-module} over $S$, if an action $(S\setminus 0)\times A\to A$ is
defined which satisfies for all $s,t\in S\setminus 0,\ a,b\in A$ the
following
conditions:
$$
s(a+b)=sa+sb,
$$
$$
st\ne 0 \Rightarrow s(ta)=(st)a.
$$

{\it A n-dimensional 0-cochain} is a partial $n$-place mapping out of
$S$ to $A$ which is defined for all $n$-tuples \mbox{$(s_1,\ldots,s_n)$,}
such that $s_1\cdot \ldots \cdot s_n\ne 0$. The coboundary operator is
given
like for the usual cohomology by the formula
\begin{eqnarray}
\partial^nf(s_1,...,s_{n+1})&=&s_1f(s_2,...,s_{n+1}) +
\displaystyle{\sum\limits_{i=1}^n}(-1)^if(s_1,...,s_is_{i+1},...,s_{n+1})
\nonumber\\
&+&(-1)^{n+1}f(s_1,...,s_n)
\nonumber
\end{eqnarray}
The equality $\partial^2=0$ is valid too; obtained cocycles (cohomology)
are called {\it  0-cocycles (0-cohomology)} and their groups are denoted by
$Z^n_0(S,A)$ (resp. $H^n_0(S,A)$).

Note that for a semigroup $T^0=T\cup0$ with an adjointed zero
$$
H^n_0(T^0,A) \cong H^n(T,A),
$$
so 0-cohomology may be considered as a generalization of Eilenberg --
MacLane cohomology.

\bigskip

Let $L$ be a finite-dimensional normal extension of a field $K$ with the
Galois group $G$, $L^{\times}$ be the multiplicative group of $L$.

By {\it a modification} $G(\star)$ of the group $G$ we mean a semigroup on
the set $G^0=G\cup 0$ with operation $\star$ such that $x\star y$ is equal
either to $xy$ or to 0, while
$$
0\star x=x\star 0=0\star 0=0
$$
and the identity of $G$ is the same for the semigroup $G(\star)$.

In other words, to obtain a modification, one must erase the contents of
some
inputs in the multiplication table of $G$ and insert there zeros so that
the
new operation would be associative.

Note some general properties of modifications. Firstly, a modification of
$G$ satisfies the weak cancellation condition: from $x\star z=y\star z\ne
0$
it follows $x=y$ and analogously for left cancellation. Secondly, let $U$
be the subgroup of invertible elements in $G(\star)$. Then its complement
$I=G(\star)\setminus U$ is a two-sided ideal. It follows from finiteness
of $G$ that $I$ is nilpotent\cite{n3}.

$L^{\times}$ is a 0-module over every modification, where elements of
the modification act on $L^{\times}$ as automorphisms of the field.
0-cohomology groups $H^2_0(G(\star),L^{\times})$ will be called
{\it components of the Brauer monoid}. In the case when operation $\star$
is defined by such a way that $x\star y=xy$ for $x,y\ne 0$, the
component of the Brauer monoid turns out the Brauer group.

Let $S=G(\star)$ and $T=G(\ast)$ be modifications of $G$. We write
$S\prec T$, if $x\ast y=0$ implies $x\star y=0$ for all $x,y\in G$.
Clearly, in this case a homomorphism is defined
$$
\varepsilon_{T,S}:H^2_0(T,L^{\times})\longrightarrow H^2_0(S,L^{\times})
$$
Since for $S\prec T\prec U$ these homomorphisms yield the equalities
$$
\varepsilon_{U,T}\varepsilon_{T,S}=\varepsilon_{U,S},
\qquad \varepsilon_{S,S}=id,
$$
one can build by the standard way \cite{cp} the semilattice of groups
$H^2_0(S,L^{\times})$ (where $S$ runs over all modifications of $G$),
which is called {\it a (relative) Brauer monoid} $M(G,L)$.

\section{Invertible elements in modifications}

As above let $S=G(\star)$ be a modification of the Galois group $G$, $U$
the
subgroup of invertible elements of $S$, $I=S\setminus U$. We shall write in
this section $xy$ instead of $x\star y$, inasmuch  as the operation of the
group $G$ will not be used. Besides let us agree to employ additive
notation
for the $G$-module $L^{\times}$.

The inclusion $U\hookrightarrow G$ induces a homomorphism
$$
\varphi: H_0^2(S,L^{\times})\rightarrow H^2(U,L^{\times})
$$
We study this homomorphism in the situation when $U$ is a normal subgroup
of $G$. Then $U$ turns out normal in $S$ too (in the meaning that $xU=Ux$
for all $x\in S$). The partition into cosets of $U$ (together with \{0\}) is
a congruence on $S$. The quotient-semigroup on this congruence will be
denoted by $S/U$.

Further, if $U\triangleleft S$ then the multiplicative group of the
subfield $P$ of all $U$-fixed elements is a $S/U$-module. The inclusion
$P^{\times}\hookrightarrow L^{\times}$ and the epimorphism $S\to S/U$
induce a homomorphism
$$
\psi: H_0^2(S/U,P^{\times})\rightarrow H_0^2(S,L^{\times})
$$

\begin{theorem} \label{2.1}
Let $U\triangleleft S$. Then the sequence
$$
0\longrightarrow H^2_0(S/U, P^{\times})
\mathop{\longrightarrow}\limits^{\psi}
H_0^2(S,L^{\times})\mathop{\longrightarrow}\limits^{\varphi}
H^2(U,L^{\times})
$$
is exact.
\end{theorem}
{\bf Proof.} We fix a system $T$ of representatives of the cosets with
respect
to $U$: $S=\{0\}\cup(\bigcup_{t\in T}Ut)$ (as usually the representative
for
$U$ is the identity). In what follows we denote by $a,b,c$ elements of $U$,
by
$x,y,z$ elements of $I\setminus 0$, by $r,s,t$ elements of $T$. If $f$ is a
(0)-cocycle then the respective element of the (0)-cohomology group is
denoted by $[f]$.

A remark must yet be made. It follows out of normality of $U$
that for every $a\in U$, $x\in I\setminus 0$ there is such $a_x\in U$ that
$ax=xa_x$. In addition $a_x$ is defined uniquely because $xa_x\ne 0$. Next,
if $a$ runs over $U$ then $a_x$ runs over it too, and  $(ab)_x=a_xb_x$ for
$a,b\in U$.

Let $[f]\in {\rm Ker}\varphi$ where $f$ is a 0-cocycle of
$Z_0^2(S,L^{\times})$. One can assume that $f(a,b)=0$ for all $a,b\in U$.
Consider the cochain of the group $U$
$$
\pi_x (a)=f(a,x)-f(x,a_x).
$$
It is a cocycle because
\begin{eqnarray}
\partial \pi_x (a,b)&=&
af(b,x)-af(x,b_x)-f(ab,x)+ f(x,a_xb_x)+f(a,x)-f(x,a_x)
\nonumber\\
&=& \partial f(a,b,x)-\partial f(a,x,b_x)-f(ax,b_x)+
f(x,a_xb_x)-f(x,a_x)
\nonumber\\
&=& -f(xa_x,b_x)+f(x,a_xb_x)-f(x,a_x)=\partial f(x,a_x,b_x)=0
\nonumber
\end{eqnarray}
Since $H^1(U,L^{\times})=0$ (see, e. g., \cite{la}), it follows from this
that
$$
\pi_x (a)=(a-1)\lambda(x),
$$
where $\lambda(x)\in L^{\times}$ for all $x\in I\setminus 0$. Set
$\lambda(a)=0$ and $g=f+\partial \lambda$. Then $g(a,b)=0$ and besides
\begin{eqnarray}
g(a,x)&=& f(a,x)-a\lambda(x)+\lambda(ax)
\nonumber\\
&=& \pi_x (a)+f(x,a-x)-a\lambda(x)+\lambda(xa_x)
\nonumber\\
&=& f(x,a-x)-\lambda(x)+\lambda(xa_x)=g(x,a_x).
\nonumber
\end{eqnarray}

Next let us set $\rho(at)=g(a,t)$ for $t\in T$ and consider the 0-cocycle
$h=g+\partial\rho$. Then
\begin{eqnarray}\label{2.eq1}
h(a,bt)&=&g(a,bt)+ag(b,t)-g(ab,t)=\partial g(a,b,t)+g(a,b)=0\\
\label{2.eq2}
h(at,b)&=&g(at,b)-g(t,a_xb)+g(a,t)=\partial g(t,a_x,b)-tg(a_x,b)=0
\end{eqnarray}
From here we obtain for $xy\ne 0$:
\begin{eqnarray}
ah(x,y)&=& h(ax,y)-h(a,xy)+h(a,x)=h(xa_x,y)
\nonumber\\
&=&xh(a_x,y)+h(x,a_xy)-h(x,a_x)=h(x,a_xy)=h(x,y(a_x)_y)
\nonumber\\
&=&-xh(y,(a_x)_y)+h(xy,(a_x)_y)+h(x,y)=h(x,y)
\nonumber
\end{eqnarray}

Hence $h(x,y)\in P^{\times}$. Besides the last calculation implies
$$
h(ax,y)=h(x,a_xy)=h(x,y)
$$
Since $a$ and $a_x$ both run over  all group $U$ we have:
$$
h(ax,y)=h(x,ay)=h(x,y)
$$
This means that $h$ defines a 0-cocycle
$\overline{h}\in Z^2_0(S/U, P^{\times})$ by the next way:
$$
\overline{h}(Us,Ut)=h(s,t)\ \  {\rm for}\  st\ne 0.
$$

So for a given 0-cocycle $f\in Z_0^2(S,L^{\times})$ by the 0-cocycles $g$
and $h$ which are cohomological to it, we construct the 0-cocycle
$\overline{h}\in Z^2_0(S/U, P^{\times})$.

We show that the correspondence $h\rightarrow \overline{h}$ extends to the
cohomology mapping. Let $h=\partial \sigma$ for some 0-cochain
$\sigma \in C^1_0(S, L^{\times})$. Since $\partial \sigma (a,b)=h(a,b)=0$,
one has $\sigma (a)=(a-1)\mu$ for the restriction of $\sigma$ on $U$, where
$\mu \in L$. Further, it follows out of (\ref{2.eq1}) and (\ref{2.eq2})
$\partial \sigma (a,x)=\partial \sigma (x,a)=0$. Therefore
\begin{eqnarray}
a[\sigma(x)-(x-1)\mu]&=&
\partial \sigma (a,x)+\sigma (ax)-\sigma (a)-xa_x\mu+a\mu
\nonumber\\
&=& \sigma (xa_x)-xa_x\mu+\mu
\nonumber\\
&=& \partial \sigma (x,a_x)-x\sigma (a_x)+\sigma (x)-xa_x\mu+\mu
\nonumber \\
&=& \sigma(x)-(x-1)\mu,
\nonumber
\end{eqnarray}
so $\sigma(x)-(x-1)\mu \in P^{\times}$. Let $\tau (g)=\sigma(g)-(g-1)\mu$
for any $g\in S\setminus 0$. Then
$\partial \tau =\partial \sigma =h$, $\tau(g)\in P^{\times}$ and in
addition
\begin{eqnarray}
\tau (ax)&=& -\partial\sigma(a,x)+a\sigma(x)+\sigma(a)-(ax-1)\mu
\nonumber\\
&=& a[\tau(x)-(x-1)\mu]+(a-1)\mu-(ax-1)\mu
\nonumber\\
&=& a\tau(x)= \tau(x),
\nonumber
\end{eqnarray}
i.\,e. 0-cochain $\tau$ is constant on cosets of $U$. Setting
$\overline{\tau}(Ut)=\tau(t)$ we get
$\overline{h}=\partial\overline{\tau}$.
Thus a homomorphism ${\rm Ker}\varphi\rightarrow H^2_0(S/U, P^{\times})$
is defined.

Now we construct an inverse map. Let $\overline{h}\in Z^2_0(S/U,
P^{\times})$.
Set $h(at,b)=h(a,bt)=0$ and $h(as,bt)=\overline{h}(Us,Ut)$
for $st\ne 0$. Then one can verify straightforward that $\partial h=0$, and
since the restriction of $h$ on $U$ equals zero, $[h]\in {\rm Ker}\varphi$.

Let $\overline{h}=\partial\overline{\gamma}$ for some 0-cochain
$\overline{\gamma}\in C^1_0(S/U, P^{\times})$.
Setting $\gamma(at)=\overline{\gamma}(Ut)$ we get $h=\partial\gamma$. Thus
we constructed the sought mapping
${\rm Ker}\varphi \to H^2_0(S/U, P^{\times})$ and proved that these groups
are isomorphic. $\blacksquare$

{\bf Remark 1.} The proved theorem generalizes the results from \cite{hls}
(where one assumed that $I^2=0$) and \cite{n3} (where one assumed that
the modification $S$ is commutative).
                                                  
{\bf Remark 2.} Indeed $S/U$ is a modification of the group $G/U$,
the Galois group of the extension $P/K$, so $H^2_0(S/U, P^{\times})$ is a
component of the respective Brauer monoid $M(G/U,P)$.

\begin{corollary} \label{2.2}
If the field $L$ is finite then
$$
H_0^2(S,L^{\times})\cong H^2_0(S/U, P^{\times})
$$
\end{corollary}

{\bf Proof.} Since in this case the group $G$ is Abelian (even cyclic)
so $U\triangleleft S$. Then $U$ is the Galois group of $L/P$ and
$H^2(U,L^{\times})$ is trivial as a Brauer group of a finite field.
$\blacksquare$

\bigskip

V.V.Kirichenko, Vozdukhoflotskij prospekt 20/1, apt.47, Kiev,
252049, Ukraine\\
e-mail: vkir@mechmat.univ.kiev.ua

B.V.Novikov, Saltovskoye shosse 258, apt.20, Kharkov, 310178, Ukraine\\
e-mail: boris.v.novikov@univer.kharkov.ua

\end{document}